\documentclass{amsart}

\newtheorem{theorem}{Theorem}[section]
\newtheorem{proposition}[theorem]{Proposition}
\newtheorem{remark}[theorem]{Remark}
\numberwithin{equation}{section}

\begin{document}

\title[Some combinatorial coincidences for standard representations]
 {Some combinatorial coincidences for standard representations of affine Lie algebras}

\author{Mirko Primc}

\address{University of Zagreb, Zagreb, Croatia}

\email{primc@math.hr}

\subjclass[2000]{17B67 (Primary); 17B69, 05A19 (Secondary).\\
\indent  Partially supported by the Croatian Science
Foundation under the Project 2634 and by the QuantiXLie Centre of Excellence, a project cofinanced
by the Croatian Government and European Union through the European Regional Development Fund---the Competitiveness and Cohesion Operational Programme (KK.01.1.1.01).}

\begin{abstract}
In this note we explain, in terms of finite dimensional representations of Lie algebras $\mathfrak{sp}_{2\ell}\subset\mathfrak{sl}_{2\ell}$,  a combinatorial coincidence of difference conditions in two constructions of combinatorial bases for standard representations of symplectic affine Lie algebras.
\end{abstract}

\maketitle
\section{Introduction}

The first of two famous Rogers-Ramanujan combinatorial identities states that the
number of partitions of $m$, written as 
\begin{equation*}
m=\sum_{j\geq 1}jf_j,
\end{equation*}
with parts $j$ congruent to $\pm 1\textrm{ mod }5$, equals the number of partitions of $m$ such that a difference between any two consecutive parts is at least two. Difference two condition can be written in terms of frequencies $f_{j}$ of parts $j$ as
\begin{equation}\label{E: difference condition 1}
f_{j}+f_{j+1}\leq 1,\qquad j\geq 1.
\end{equation}
J. Lepowsky and S. Milne discovered in \cite{LM} that the generating function of partitions satisfying congruence condition is, up to a certain ``fudge'' factor $F$, the principally specialized character of certain level $3$ standard module $L(\Lambda)$ for affine
Lie algebra $\widehat{\mathfrak{sl}}_2$.  Lepowsky and R.L. Wilson
realized that the factor $F$ is a character of the Fock space for the principal Heisenberg subalgebra of  $\widehat{\mathfrak{sl}}_2$,
and that the generating function of partitions satisfying difference two condition is the principally specialized character of the vacuum space of $L(\Lambda)$ for the action of principal Heisenberg subalgebra.
In a series of papers---see \cite{LW} and the references therein---Lepowsky and Wilson  gave a Lie-theoretic proof of both Rogers-Ramanujan identities by constructing combinatorial  bases of vacuum spaces for the principal Heisenberg subalgebra
parametrized by partitions satisfying difference two conditions. Very roughly speaking, in the case of the first Rogers-Ramanujan identity, the vacuum
space $\Omega$ is spanned by monomial vectors of the form
\begin{equation}\label{E: LW argument 1}
Z(-s)\sp{f_s}\dots Z(-2)\sp{f_2}Z(-1)\sp{f_1}v_\Lambda,\quad s\geq 0, \ f_j\geq 0,
\end{equation}
where $v_\Lambda\in\Omega$ is a highest weight vector of $L(\Lambda)$ and $Z(j)$ are certain $\mathcal Z$-operators. The degree of monomial vector
(\ref{E: LW argument 1}) is
\begin{equation*}
-m=-\sum_{j=1}\sp s jf_j.
\end{equation*}
Lepowsky and Wilson discovered that $\mathcal Z$-operators $Z(j)$ on $\Omega$ satisfy certain relations, roughly of the form
\begin{equation}\label{E: LW relations}
\begin{aligned}
&Z(-j)Z(-j)+2\sum_{i>0}Z(-j-i)Z(-j+i)\approx 0,\\
&Z(-j-1)Z(-j)+\sum_{i>0}Z(j-1-i)Z(-j+i)\approx 0,
\end{aligned}
\end{equation}
so we may replace {\it the leading terms}
\begin{equation}\label{E: LW argument 3}
Z(-j)Z(-j) \ \text{and} \  Z(-j-1)Z(-j)
\end{equation}
of relations (\ref{E: LW relations}) with ``higher terms'' $Z(-j-i)Z(-j+i)$ and $Z(-j-1-i)Z(-j+i)$, $i>0$, and reduce the spanning set
(\ref{E: LW argument 1}) of $\Omega$ to a spanning set
\begin{equation}\label{E: LW argument 4}
Z(-s)\sp{f_s}\dots Z(-2)\sp{f_2}Z(-1)\sp{f_1}v_0,\quad s\geq 0,\  f_{j}+f_{j+1}\leq 1\ \text{for all} \  j\geq 1.
\end{equation}
By invoking the product formula for principally specialized character of $\Omega$ and the first Rogers-Ramanujan identity, we see that vectors
in the spanning set (\ref{E: LW argument 4}) are in fact a basis of $\Omega$. Lepowsky and Wilson proved directly linear independence of the spanning set (\ref{E: LW argument 4}) and, by reversing the above argument, proved the combinatorial identity. 

Lepowsky-Wilson's approach has been applied in many different situations, here we are interested in particular two constructions of combinatorial bases of level $k\geq1$ standard modules for affine Lie algebras of higher ranks in ``homogeneous picture''---both constructions start with a relation of the form
\begin{equation}\label{E: MP argument 3}
\sum_{j_1+\cdots +j_{k+1}=m}x(j_1)\cdots
x(j_{k+1})=0,
\end{equation}
(see (\ref{E: defining relation}) below) with the leading term
\begin{equation}\label{E: MP argument 4}
 x(-j-1)^b  x(-j)^a,
\end{equation}
$a+b=k+1$ and $(-j-1)b+(-j)a=m$. For $k=1$ this is analogous to (\ref{E: LW relations}) and (\ref{E: LW argument 3}), and again we may replace the leading terms (\ref{E: MP argument 4}) with higher terms in
(\ref{E: MP argument 3}) and reduce a Poincar\' e-Birkhoff-Witt spanning set of monomial vectors in the given representation space
to a smaller spanning set. However, this spanning set is not a basis and,
in the case we consider, all the relations needed to reduce the Poincar\' e-Birkhoff-Witt spanning set to a basis are obtained from (\ref{E: MP argument 3}) by
the adjoint action of the finite dimensional Lie algebra $\mathfrak{sp}_{2\ell}$ or $\mathfrak{sl}_{2\ell}$, and all the leading terms which we replace by higher terms are obtained by the adjoint action of the same Lie algebra on the leading terms (\ref{E: MP argument 4}). In analogy to the Rogers-Ramanujan case above, we say that the monomial vectors in constructed combinatorial basis {\it satisfy difference conditions for a given module}.

In this note we explain, in terms of finite dimensional representations of Lie algebras $\mathfrak{sp}_{2\ell}\subset\mathfrak{sl}_{2\ell}$,  a combinatorial coincidence of difference conditions in two constructions \cite{BPT, T3} and \cite{PS1} of combinatorial bases for standard representations of symplectic affine Lie algebra. In particular, Proposition \ref{proposition} below and Propositions 1 and 2 in \cite{BPT} provide an alternative proof of Theorem 6.1 in \cite{PS2}.

\section{Affine Lie algebras}

We are interested mainly in symplectic affine Lie algebra, but it will be convenient to have a bit more general notation for an  affine Lie algebra $\tilde{\mathfrak g}$.

\subsection{Affine Lie algebras}
Let ${\mathfrak g}$ be a simple complex Lie algebra, $\mathfrak h$
a Cartan subalgebra, $R$ the corresponding root system and $\theta$  the maximal root with respect to some fixed basis of the root system. For each root $\alpha$ we fix a root vector $x_\alpha$ in $\mathfrak g$. Via a symmetric invariant bilinear form $\langle \ , \ \rangle$ on ${\mathfrak g}$ we identify $\mathfrak h$ and $\mathfrak h^*$ and we assume
that $\langle \theta , \theta \rangle=2$.  Set
$$
\hat{\mathfrak g} =\coprod_{j\in\mathbb Z}{\mathfrak g}\otimes
t^{j}+\mathbb C c, \qquad \tilde{\mathfrak g}=\hat{\mathfrak
g}+\mathbb C d.
$$
Then $\tilde{\mathfrak g}$ is the associated untwisted affine
Kac-Moody Lie algebra (cf. \cite{K}) with the commutator
$$
[x(i),y(j)]=[x,y](i+j)+i\delta_{i+j,0}\langle x,y\rangle c.
$$
Here, as usual, $x(i)=x\otimes t^{i}$ for $x\in{\mathfrak g}$ and
$i\in\mathbb Z$, $c$ is the canonical central element, and
$[d,x(i)]=ix(i)$.  We identify ${\mathfrak g}$ and
${\mathfrak g}\otimes 1$.
Let $B$ be a basis of $\mathfrak{g}$ consisting of root vectors and elements of $\mathfrak{h}$. Set
$$
\hat B=\{x(n)\mid x\in B, n\in\mathbb Z\}.
$$
Then $\hat B\cup\{c\}$ is a basis of $\hat{\mathfrak{g}}$.
A given linear order $\preceq$ on $B$ we extend to a linear order on $\hat B$ by
$$
x(n)\prec y(m)\quad\text{iff}\quad n<m \ \text{or} \  n=m \ \text{and} \   x\prec y.
$$

\subsection{$\mathbb Z$-gradings of $\mathfrak g$ by minuscule coweights}
Let $R\sp\vee$ be the dual root system of $R$ and fix a minuscule coweight $\omega \in P(R^\vee)$. Set $\Gamma = \{\,\alpha \in R \mid \omega(\alpha\,) = 1\}$ and
$$
{\mathfrak g}_0 = {\mathfrak h} + \sum_{\omega(\alpha)=0}\,
{\mathfrak g}_\alpha, 
\qquad  {\mathfrak g}_{\pm1} = \sum_{\alpha \in \pm \Gamma}\, {\mathfrak
g}_\alpha,
$$
$$
 \hat{\mathfrak{g}}_0=\mathfrak{g}_0 \otimes \mathbb{C}[t,t^{-1}]\oplus \mathbb{C} c , 
\qquad
\hat{\mathfrak{g}}_{\pm1}=\mathfrak{g}_{\pm1}\otimes\mathbb{C}[t,t^{-1}].
$$
Then on $\mathfrak g$ and $\hat{\mathfrak{g}}$ we have $\mathbb{Z}$-gradings
$$
\mathfrak{g} = \mathfrak{g}_{-1} + \mathfrak{g}_0 + \mathfrak{g}_1\quad\text{and}\quad \hat{\mathfrak{g}} = \hat{\mathfrak{g}}_{-1} + \hat{\mathfrak{g}}_0 + \hat{\mathfrak{g}}_1 .
$$ 
Let ${\mathfrak h'}_0\subset{\mathfrak h}$ be a Cartan subalgebra of ${{\mathfrak g}_0}'=[{\mathfrak g}_0,{\mathfrak g}_0]$. Note that  ${\mathfrak g}_0 ={{\mathfrak g}_0}' +\mathbb C\omega$ is a reductive Lie algebra,  ${\mathfrak g}_{\pm1}$ are commutative subalgebras of $\mathfrak g$ and $\mathfrak g _0$-modules, and that
$\hat{\mathfrak{g}}_{\pm1}$ are commutative subalgebras of $\hat{\mathfrak{g}}$ and  $\hat{\mathfrak{g}}_0$-modules.

For a classical $\mathfrak g$ all possible minuscule coweights $\omega$ are listed on Dynkin diagrams below: the black dot on a diagram denotes the simple root $\alpha$ for which $\omega(\alpha)=1$ and the rest is the Dynkin diagram of the semisimple Lie algebra ${{\mathfrak g}_0}' $.
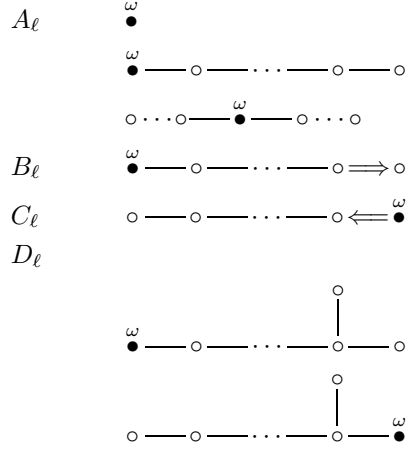
\smallskip

\def\hcrta{\operatorname{\vrule width3.5ex height.6ex depth-.5ex}}
\def\vcrtao{\vbox{\hbox{$\circ$}\kern-1ex 
\hbox{\kern.5ex\vrule width.1ex height2.2ex depth1ex}\kern-1ex\hbox{$\circ$}}}
$$
\begin{aligned}
& A_\ell\quad &&&&\,\overset{\omega}{\bullet} && \\
&  &&
&&\overset{\omega}{\bullet}\hcrta\circ\hcrta\cdots\hcrta\circ\hcrta\circ
\qquad && \\
& &&
&&\,{\circ} 
\cdots\!\circ \!\hcrta \overset{\omega}{\bullet} \hcrta
\circ
\cdots \circ
\qquad &&  \\
&B_\ell&&
&&\overset{\omega}{\bullet}\hcrta\circ\hcrta
\cdots\hcrta\circ\!\Longrightarrow\!\circ
\qquad &&  \\
& C_\ell &&
&&\circ\hcrta\circ\hcrta
\cdots\hcrta\circ\!\Longleftarrow\!\overset{\omega}{\bullet}
\quad && \\
& D_\ell&&  && && \\
& && 
&&\overset{\omega}{\bullet}\hcrta\circ\hcrta \cdots\hcrta 
\vcrtao\hcrta\circ
&&\\
&&&        &&\circ\hcrta\circ\hcrta\cdots\hcrta 
\vcrtao\hcrta
\overset{\omega}{\bullet}
&&  
\end{aligned}
$$

\section{Standard modules and Feigin-Stoyanovsky's type subspaces}

In this section we describe monomial bases of standard $\tilde{\mathfrak{g}}$-modules $L(k\Lambda_0)$ and monomial bases of Feigin-Stoyanovsky's type subspaces $W(k\Lambda_0)$.

\subsection{Standard modules and Feigin-Stoyanovsky's type subspaces}
As usual, we denote by $\Lambda_0, \dots, \Lambda_{\ell}$ the fundamental weights of $\tilde{\mathfrak{g}}$. For a given integral dominant weight $\Lambda = k_0\Lambda_0 + \cdots + k_{\ell} \Lambda_{\ell}$ we denote by $L(\Lambda)$ the standard (i.e. integral highest weight) $\tilde{\mathfrak{g}}$-module with the highest weight $\Lambda$, by $v_{\Lambda}$ a fixed highest weight vector and by  $k=\Lambda(c)$ the level. 
\smallskip

For a standard $\tilde{\mathfrak{g}}$-module $L(\Lambda)$ and a given $\mathbb{Z}$-grading on $\hat{\mathfrak{g}}$ as above, we define the corresponding Feigin-Stoyanovsky's type subspace $W(\Lambda)$ as 
$$
 W(\Lambda) = U(\hat{\mathfrak{g}}_1)  v_{\Lambda}\subset
L(\Lambda).
 $$
Note that $U(\hat{\mathfrak{g}}_1) \cong S(\hat{\mathfrak{g}}_1) $ since $\hat{\mathfrak{g}}_1$ is commutative.

\subsection{Monomials in elements of $\hat B$ and colored partitions} 
Let $B$ be a basis of $\mathfrak{g}$ consisting of root vectors and elements of $\mathfrak{h}$, with a linear order $\preceq$. 
We shall write ordered product monomials in 
$\hat B=\{x(n)\mid x\in B, n\in\mathbb Z\}$ with factors $x_i(n_i)$  in ascending order, i.e. as
\begin{equation}\label{E: ordered product monomial}
\prod_{i=1}\sp rx_i(n_i)= x_1(n_1)x_2(n_2)\dots x_r(n_r),\quad x_1(n_1)\preceq x_2(n_2)\preceq\dots \preceq x_r(n_r).
\end{equation}
Sometimes it is convenient to think of the ordered sequence
\begin{equation}\label{E: colored partitions}
\pi=(x_1(n_1)\preceq x_2(n_2)\preceq\dots \preceq x_r(n_r))
\end{equation}
as a colored partition, i.e., as a ``plain partition'' $n_1\leq n_2\leq\dots\leq n_r$ colored with ``colors'' $x_1,  x_2, \dots , x_r$. We shall say that $r$ is the length of $\pi$ and that  $n_1+n_2+\dots+n_r$ is the degree of $\pi$.
Then the monomial in (\ref{E: ordered product monomial}) is the corresponding ordered product $u(\pi)$ in the enveloping algebra
\begin{equation}\label{E: monomial vector u(pi)}
u(\pi)= x_1(n_1)x_2(n_2)\dots x_r(n_r)\in  U(\hat{\mathfrak{g}}).
\end{equation}
It will is also be convenient to think of the colored partition $\pi$ as a monomial in a symmetric algebra of $\hat{\mathfrak{g}}$, i.e.
$$
\pi= x_1(n_1)x_2(n_2)\dots x_r(n_r)\in  S(\hat{\mathfrak{g}}).
$$
We denote with $\mathcal P$ the set of all colored partitions; it is a monoid with the unit $1$---the colored partition with no parts (of length $0$ and degree $0$). The submonoid $\mathcal P_{<0}$ is the set of all colored partitions (\ref{E: colored partitions}) such that $n_1\leq n_2\leq\dots\leq n_r<0$.

When $\omega$ is a minuscule coweight and $B_1$ an ordered basis of $\mathfrak g_1$, we extend the order to a basis ${\hat B}_1=\{x_\gamma(n)\mid \gamma\in \Gamma, n\in\mathbb Z\}$ of $\hat{\mathfrak g}_1$ and, in a similar way as above, we consider colored partitions and ordered monomials in $U(\hat{\mathfrak g}_1)\cong S(\hat{\mathfrak g}_1))$. We denote the corresponding monoids as $\mathcal P\sp1$ and $\mathcal P\sp1_{<0}$.

\subsection{Well order on colored partitions}
We extend the order $\preceq $ on $\hat B$ to the order on the set of colored partitions (\ref{E: colored partitions}) first by comparing the lengths---shorter partitions are higher; then by comparing the degrees---partitions with greater degree are higher; then by comparing ``plain partitions'' $n_1\leq n_2\leq\dots\leq n_r$ in the reverse lexicographical order and, finally,  by comparing ``the colorings'' $x_1,  x_2, \dots , x_r$ in the reverse lexicographical order. 
This order has several good properties:

\noindent $\bullet$ On the set of partitions with bounded length and fixed degree $\preceq $ is a well order, so in our arguments we can use induction.

\noindent $\bullet$ For $x_1,x_2,\dots, x_r\in \hat B$ and any permutation $\sigma$ two monomials $x_1x_2\dots x_r$ and $x_{\sigma(1)}x_{\sigma(2)}\dots x_{\sigma(r)}$ in $U( \hat{\mathfrak{g}})$ differ by a linear combination of $u(\pi)$ with length of $\pi$ strictly less then $r$. So in inductive arguments ordered monomials ``behave as commutative monomials''.

\noindent $\bullet$ \  $\kappa\preceq \lambda$ implies $\kappa\pi\preceq \lambda\pi$ for colored partitions $\kappa, \lambda, \pi$ in  $S(\hat{\mathfrak{g}})$. This property enables us to replace the leading terms of  relations with higher terms and use induction.

\subsection{Relations on  $L(\Lambda)$ and $W(\Lambda)$}
 On level $k$ standard module $L(\Lambda)$
we have vertex operator relations (cf. \cite{LP, MP2, P1})
\begin{equation}\label{E: defining relation}
x_\theta(z)^{k+1}= \sum_{n\in\mathbb Z}\Big(\sum_{j_1+\dots+j_{k+1}=n}
x_\theta(j_1)\dots x_\theta(j_{k+1})\Big)z^{-n-k-1}=0.
\end{equation}
With the adjoint action of $\mathfrak g$ we get a finite dimensional $\mathfrak g$-module
$$
U(\mathfrak g)\cdot x_\theta(z)^{k+1}\cong L_{\mathfrak g}((k+1)\theta).
$$
The set $\bar{R}_k$ of all coefficients in $U(\mathfrak g)\cdot x_\theta(z)^{k+1}$ we call {\it relations on  $L(\Lambda)$}. 
\medskip

For a minuscule coweight $\omega$ we have $\omega(\theta)=1$. Hence $x_\theta\in \mathfrak g_1$ and the vertex operator relation (\ref{E: defining relation}) is a relation on the Feigin-Stoyanovsky's type subspace $W(\Lambda)$. Since $\mathfrak g_1$ is a $\mathfrak g_0$-module, by the adjoint action of $\mathfrak g_0$ on the vertex operator relation (\ref{E: defining relation}) we get a finite dimensional $\mathfrak g_0$-module
$$
U(\mathfrak g_0)\cdot x_\theta(z)^{k+1}\cong L_{{\mathfrak g_0}'}((k+1)\theta\vert_{\mathfrak h'_0} ),
$$
where $\mathfrak h'_0\subset\mathfrak h$ is a Cartan subalgebra of ${\mathfrak g_0}'$. The set $\bar{R}\sp0_k\subset \bar{R}_k$ of all coefficients in $U(\mathfrak g_0)\cdot x_\theta(z)^{k+1}$ we call {\it relations on  $W(\Lambda)$}.

\subsection{The leading terms of relations on  $L(\Lambda)$ and $W(\Lambda)$}
Every relation $r\in\bar R$, $r\neq0$, can be written as
$$
r=c_\rho u(\rho)+\sum_{\pi\succ\rho}c_\pi u(\pi),\qquad c_\rho\neq 0.
$$
We say that the colored partition $\rho$ is {\it the leading term of $r$} and we write $\rho=\ell t(r)$, and sometimes we shall also say that the ordered monomial $u(\rho)$ is the leading term of the relation $r=0$. 
The leading terms of coefficients of (\ref{E: defining relation}) are colored partitions of the form
\begin{equation}\label{E: leading terms of r(k+1)theta}
 x_\theta(-j-1)^b  x_\theta(-j)^a\in S(\hat{\mathfrak{g}}),
\end{equation}
$a+b=k+1$ and $(-j-1)b+(-j)a=-n$. All leading terms of relations $\bar{R}_k$ and $\bar{R}\sp0_k$ are obtained as leading terms of finite dimensional spaces
\begin{equation}\label{E: all leading terms}
\ell t\left(U(\mathfrak g)\cdot x_\theta(-j-1)^b  x_\theta(-j)^a\right)\quad\text{and}\quad
\ell t\left(U(\mathfrak g_0)\cdot x_\theta(-j-1)^b  x_\theta(-j)^a\right).
\end{equation}
 
\subsection{Monomial bases of  Feigin-Stoyanovsky's type subspaces $W(k\Lambda_0)$} 
By Poincar\' e-Birkhoff-Witt theorem the set of monomial vectors
$$
u(\pi)v_{k\Lambda_0},\qquad \pi\in \mathcal P\sp1_{<0}
$$
is a spanning set for  Feigin-Stoyanovsky's type subspace $W(k\Lambda_0)$. By expressing the leading terms of relations in $\bar{R}\sp0_k$ with higher monomials, we can reduce this PBW spanning set of $W(k\Lambda_0)$ to a spanning set
\begin{equation}\label{E: basis of FS subspace}
u(\pi)v_{k\Lambda_0},\qquad \pi\in \mathcal P\sp1_{<0}, \quad \pi\not\in \mathcal P\sp1
\ell t\left( \bar{R}\sp0_k \right).
\end{equation}
This set is a basis in many cases with a suitable choice of order $\preceq$ and $B$: for all positive integer levels $k$ and $\mathfrak g$ of the type $A_\ell$ for all $\ell$ and all minuscule coweights (see \cite{CLM1, CLM2, FJLMM, FS, P1, P3, T1, T2}) and for $\mathfrak g$ of the type $C_\ell$ for all $\ell$ (see  \cite{BPT, P4}); for level $k=1$ for all classical $\mathfrak g$  for all $\ell$ and all minuscule coweights (see \cite{P3}); and for $D_4$ level $k=1$ and $2$  (see \cite{B}). The monomial bases of  Feigin-Stoyanovsky's type subspaces are used in the construction of semi-infinite monomial bases of entire standard modules (cf. \cite{T3}) .

Different methods are used to prove linear independence of the spanning set (\ref{E: basis of FS subspace}), but in all proofs the explicit combinatorial description of colored partitions $\pi$ which satisfy the difference condition $ \pi\not\in \mathcal P\sp1
\ell t\left( \bar{R}\sp0_k \right)$ is needed. For $A_1$ and $k\geq1$ partitions $\pi$ have only one color and the difference conditions are the difference two conditions which appear in Rogers-Ramanujan and Gordon combinatorial identities. For $A_\ell$ and $\omega=\omega_1$ colored partitions $\pi$ appear as $(k,\ell+1)$-admissible configurations. For all classical Lie algebras $\mathfrak g$ and all minuscule coweights difference conditions for $W(\Lambda_0)$ are given by the energy function of a perfect crystal (cf. \cite{KKMMNN}) corresponding to the $\mathfrak g_0$-module $\mathfrak g_1$.

\subsection{Monomial bases of standard modules} By Poincar\' e-Birkhoff-Witt theorem the set of monomial vectors
$$
u(\pi)v_{k\Lambda_0},\qquad \pi\in \mathcal P_{<0}
$$
is a spanning set for the standard module $L(k\Lambda_0)$. By expressing the leading terms of relations in $\bar{R}_k$ with higher monomials, we can reduce this PBW spanning set of $L(k\Lambda_0)$ to a spanning set
\begin{equation}\label{E: basis of standard module}
u(\pi)v_{k\Lambda_0},\qquad \pi\in \mathcal P_{<0}, \quad \pi\not\in \mathcal P
\ell t\left( \bar{R}_k \right).
\end{equation}
This set is a basis in several cases with a suitable choice of order $\preceq$ and $B$: for all positive integer levels $k$ and $\mathfrak g$ of the type $A_1$  (see \cite{FKLMM, F, MP1, MP2}) and for level $k=1$ and $\mathfrak g$ of the type  $C_\ell$ for all $\ell\geq 2$
(see \cite{PS1, S}). In this setting we can think $A_1\cong C_1$, so (\ref{E: basis of standard module}) is a basis of $L(k\Lambda_0)$ for $\mathfrak g$  of the type  $C_\ell$ for $k=1$ or $\ell=1$.
In \cite{PS2} we conjecture that  (\ref{E: basis of standard module}) is a basis for all $k\geq1$ and $\ell\geq1$.

\section{Some combinatorial coincidences for two constructions}

In this section we describe the combinatorial coincidence of the leading terms $\ell t\left( \bar{R}\sp0_k \right)$ for $C_{2\ell}$ and the leading terms $\ell t\left( \bar{R}_k \right)$ for $C_{\ell}$ which happens for all $k\geq 1$ and $\ell\geq 1$.

\subsection{Combinatorial coincidence of  $\ell t\left( \bar{R}\sp0_k \right)$ for $B_{2}$ and $\ell t\left( \bar{R}_k \right)$ for $A_1$}
Let
$\mathfrak g$ of the type $B_2\cong C_2$. Then ${\mathfrak g_0}'$ is of the type $A_1$. We identify the root system with
$$R=\{\pm(\varepsilon_1 - \varepsilon_2), \pm(\varepsilon_1 + \varepsilon_2), \pm\varepsilon_1,  \pm\varepsilon_2)\}.
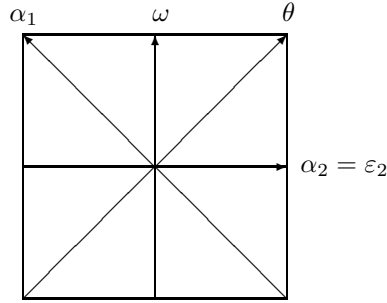
$$ 
\smallskip

\begin{center} 
\begin{picture}(100,100)
\put (0,0) {\line(1,0){100}}
\put (0,0) {\vector (1,1){100}}
\put (50,0) {\vector(0,1){100}}
\put (100,0) {\line(0,1){100}}
\put (100,0) {\vector(-1,1){100}}
\put (0,50) {\vector(1,0){100}}
\put (105,48){$\alpha_{2}=\varepsilon_2$} 
\put (0, 100) {\line(0,-1){100}}
\put (100,100) {\line (-1,0){100}}
\put (98,105){$\theta$}
\put (49,105){$\omega$}
\put (-5,105){$\alpha_{1}$}
\end{picture}
\end{center}

\noindent
Then $\omega=\omega_1$ and \ $\Gamma=\{\varepsilon_1 - \varepsilon_2, \varepsilon_1, \varepsilon_1+\varepsilon_2\}$ with the the corresponding root vectors in $\mathfrak g_1$ denoted as  $x_{\underline{2}}\,, x_0\,, x_2$ with the order $x_{\underline{2}}\prec x_0\prec  x_2$. The consecutive action of $x_{-\varepsilon_2}$ in $\mathfrak g_0$ on (\ref{E: leading terms of r(k+1)theta}) gives the leading terms

\begin{equation}\label{E: leading terms FS B2}
\begin{aligned}
&x_{{2}}(-j-1)\sp{a_{j+1}}x_0(-j)\sp{b_{j}}x_2(-j)\sp{a_j}\,,
\quad \ \,  \qquad a_{j+1}+b_{j}+a_{j}= k+1,\\
&x_{{2}}(-j-1)\sp{a_{j+1}}x_{\underline{2}}(-j)\sp{c_{j}}x_0(-j)\sp{b_j}\,,
\quad \ \, \qquad a_{j+1}+c_{j}+b_{j}= k+1,\\
&x_{{0}}(-j-1)\sp{b_{j+1}}x_2(-j-1)\sp{a_{j+1}}x_{\underline{2}}(-j)\sp{c_j}\,,
\quad b_{j+1}+a_{j+1}+c_{j}= k+1,\\
&x_{\underline{2}}(-j-1)\sp{c_{j+1}}x_0(-j-1)\sp{b_{j+1}}x_{\underline{2}}(-j)\sp{c_j}\,,
\quad c_{j+1}+b_{j+1}+c_{j}= k+1,
\end{aligned}
\end{equation}
and the spanning set 
 \begin{equation}\label{E: spanning set for FS B2}
\prod_{j<0}x_{\underline{2}}(-j)\sp{c_j}x_0(-j)\sp{b_j}x_2(-j)\sp{a_j}v_{k\Lambda_0}
\end{equation}
of monomial vectors in $W_{B^{(1)}_2}(k\Lambda_0)$ is reduced to a basis (\ref{E: basis of FS subspace}) if we impose the difference conditions
\begin{equation}\label{E: B2 FS difference conditions}
\begin{aligned}
c_{j+1}+b_{j+1}+c_{j}&\leq k,\\
b_{j+1}+a_{j+1}+c_{j}&\leq k,\\
a_{j+1}+c_{j}+b_{j}&\leq k,\\
a_{j+1}+b_{j}+a_{j}&\leq k
\end{aligned}
\end{equation}
for all $j\in\mathbb Z$ (see \cite{BPT, P4}). It is clear that the same argument for leading terms applies in the case of the standard module $L_{A^{(1)}_1}(k\Lambda_0)$ (see \cite{FKLMM, MP2}) if we use the standard ordered basis $x_{-\alpha}\prec h\prec  x_\alpha$ of $\mathfrak{sl}_2$ and identify the corresponding bases elements
\begin{equation}\label{E: identification of bases for B2 and A1}
\begin{aligned}
&x_{-\alpha}\leftrightarrow x_{\varepsilon_1 - \varepsilon_2}, \  h\leftrightarrow x_{\varepsilon_1}, \  x_{\alpha}\leftrightarrow x_{\varepsilon_1 + \varepsilon_2}\qquad\text{and},\quad\text{for all} \ j\in\mathbb Z,\\
&x_{-\alpha}(j)\leftrightarrow x_{\varepsilon_1 - \varepsilon_2}(j), \  h(j)\leftrightarrow x_{\varepsilon_1}(j), \  x_{\alpha}(j)\leftrightarrow x_{\varepsilon_1 + \varepsilon_2}(j).\\
\end{aligned}
\end{equation}

{\it Therefore, the Feigin-Stoyanovsky's type subspace $W_{B^{(1)}_2}(k\Lambda_0)$  and the standard module $L_{A^{(1)}_1}(k\Lambda_0)$ have the same combinatorial description of the  leading terms of relations, and the same combinatorial description of the difference conditions for colored partitions parameterizing monomial bases.}

\subsection{Simple Lie algebra $\mathfrak g$ of type $C_\ell$}

We fix a simple Lie algebra $\mathfrak{g}$ of the type $C_\ell$, $\ell\geq 2$. For a given Cartan subalgebra $\mathfrak h$ and the corresponding
root system $R$ we can write
\begin{equation*}
R = \{\pm(\varepsilon_i\pm\varepsilon_j) \mid i,j=1,...,\ell\}\setminus\{0\}.
\end{equation*}
We choose simple roots as in \cite{Bo}
\begin{equation*}
\alpha_1= \varepsilon_1-\varepsilon_2,  \ \alpha_2=\varepsilon_2-\varepsilon_3, \  \cdots \ \alpha_{\ell-1}=\varepsilon_{\ell-1}-\varepsilon_{\ell},
 \  \alpha_\ell = 2\varepsilon_\ell.
\end{equation*}
Then $\theta=2 \varepsilon_1$. For each root $\alpha$ we choose a root vector $X_{\alpha}$ such that $[X_{\alpha},X_{-\alpha}]=\alpha^{\vee}$. For the root vectors
$X_{\alpha}$ we shall use the following notation:
$$\begin{array}{ccc}
X_{ij}\quad \text{or just}\quad ij &  \text{if}\   &\alpha =\varepsilon_i + \varepsilon_j\ , \ i\leq j\,,\\
X_{\underline{i}\underline{j}}\quad \text{or just}\quad \underline{i}\underline{j} & \ \text{if}\  &\alpha =-\varepsilon_i - \varepsilon_j\ , \ i\geq j\,,\\
X_{i\underline{j}}\quad \text{or just}\quad i \underline{j} & \ \text{if}\  &\alpha =\varepsilon_i - \varepsilon_j\ , \ i\neq j\,.\\
\end{array}
$$
With the previous notation $x_\theta=X_{11}$. We also write for $i=1, \dots, \ell$
$$
X_{i\underline{i}}=\alpha_i^{\vee}\ \text{or just}\ i\underline{i} \,.
$$
These vectors $X_{ab}$ form a basis $B$ of $\mathfrak g$ which we shall write in a triangular scheme. For example, for $\ell=3$ the basis $B$ is
$$\begin{array}{cccccc}
11 &  &&  & & \\
12 & 22 & & & & \\
13 & 23 & 33 & & & \\
1\underline{3} & 2\underline{3} & 3\underline{3} & \underline{3}\underline{3}  & & \\
1\underline{2} & 2\underline{2} & 3\underline{2} & \underline{3}\underline{2}  & \underline{2}\underline{2}& \\
1\underline{1} & 2\underline{1} & 3\underline{1} & \underline{3}\underline{1}  & \underline{2}\underline{1} & \underline{1}\underline{1}.
\end{array}
$$
In general for the set of indices $\{1,2,\cdots ,\ell,\underline{\ell},\cdots ,\underline{2},\underline{1}\}$ we use order
\begin{equation*}
1\succ 2\succ\cdots\succ \ell-1\succ \ell\succ \underline{\ell} \succ  \underline{\ell-1} \succ\cdots \succ \underline{2} \succ \underline{1}
\end{equation*}
and a basis element $X_{ab}$ we write in $a^{th}$ column and $b^{th}$ row,
\begin{equation}\label{E: basis of g}
B=\{X_{ab}\mid b\in\{1,2,\cdots ,\ell,\underline{\ell},\cdots ,\underline{2},\underline{1}\},\ a\in\{1,\cdots ,b\}\} .
\end{equation}
By using (\ref{E: basis of g}) we define on the basis $B$ the corresponding  lexicographical order
\begin{equation*}
X_{ab}\succ X_{a'b'} \ if\ a\succ a'\ or \ a=a' \ and\ b\succ b'\ .
\end{equation*}
In other words, $X_{ab}$ is larger than $X_{a' b'}$ if $X_{a' b'}$ lies in a column $a'$ to the right of the column $a$, or $X_{ab}$ and
$X_{a' b'}$ are in the same column $a=a'$, but $X_{a'b'}$  is below $X_{ab}$.
For $\ell=1$ we get $C_1\cong A_1$ with the basis $B$
$$\begin{array}{cc}
11 &   \\
1\underline{1} & \underline{1}\underline{1} . 
\end{array}
$$
\smallskip

For a simple Lie algebra $\mathfrak{g}$ of type $C_\ell$, $\ell\geq 2$, we have the minuscule weight 
$$
\omega=\omega_\ell=\tfrac12(\varepsilon_1+\dots+\varepsilon_{\ell}),
$$
the corresponding 
 $\Gamma=\{\varepsilon_i + \varepsilon_j\mid 1\leq i\leq j\leq \ell\}$ and the basis
 \begin{equation}\label{E: basis of g1}
B_1=\{X_{\gamma}\mid \gamma\in\Gamma\}=\{X_{ij}\mid 1\leq i\leq j\leq \ell\}\subset B 
\end{equation}
 of $\mathfrak g_1$. For $C_3$ the basis $B_1$ can be written as  as the ``upper'' triangle in B:
$$\begin{array}{ccc}
11 &  & \\
12 & 22 &  \\
13 & 23 & 33.  
\end{array}
$$
In the case $C_1$ the ``upper'' triangle in $B$ is $B_1=\{11\}$.
In the case $C_{\ell}$, $\ell\geq2$, the subalgebra ${\mathfrak g_0}'$ is a simple Lie algebra of the type $A_{\ell-1}$ with a Cartan subalgebra $\mathfrak h'_0\subset \mathfrak h_0$ and the corresponding root system 
$$
\{\pm(\varepsilon_i - \varepsilon_j)\mid 1\leq i<j\leq \ell\}
$$
with the simple roots 
\begin{equation*}
\alpha_1= \varepsilon_1-\varepsilon_2,  \ \alpha_2=\varepsilon_2-\varepsilon_3, \  \cdots \ \alpha_{\ell-1}=\varepsilon_{\ell-1}-\varepsilon_{\ell}.
\end{equation*}

\subsection{The leading terms of relations for $W_{C^{(1)}_{\ell}}(k\Lambda_0)$}
For the finite dimensional ${{\mathfrak g}_0}'$-module of relations we have
\begin{equation}\label{E: relations for FS Cell}
U({\mathfrak g}_0)\cdot x_\theta(z)\sp{k+1}\cong L_{A_{\ell-1}}((k+1)\theta\vert_{\mathfrak h'_0})=L_{A_{\ell-1}}(2(k+1)\omega_1).
\end{equation}
Note that ${\mathfrak g_0}'$-module $L_{A_{\ell-1}}(\omega_1)$ is the vector representation  $\mathbb C\sp\ell$  of $\mathfrak{sl}_\ell$ with the canonical basis $e_1, \dots, e_\ell$, the module $L_{A_{\ell-1}}(2\omega_1)$ is isomorphic to the symmetric power $S\sp2(\mathbb C\sp\ell)$  with a basis 
$e_ie_j\leftrightarrow X_{ij}$\,, and the module (\ref{E: relations for FS Cell}), as a vector space, is isomorphic to $S\sp{2(k+1)}(\mathbb C\sp\ell)$ with the  canonical basis $e_1\sp{m_{1}}, \dots, e_\ell\sp{m_{\ell}}$, $m_1+\dots+m_{\ell}=2(k+1)$, which we view as multisets $\{1\sp{m_{1}}, \dots, \ell\sp{m_{\ell}}\}$.
By using the action of the root vectors $x_{-\alpha_1},  \dots, x_{-\alpha_{\ell-1}}$ in $\mathfrak g_0$ on (\ref{E: defining relation}) for every multiset $\{1\sp{m_1}, \dots, \ell\sp{m_\ell}\}$, $m_1+\dots+m_{\ell}=2(k+1)$, we obtain a relation for Feigin-Stoyanovsky's type subspace $W_{C^{(1)}_{\ell}}(k\Lambda_0)$ (see \cite{BPT}, Proposition 1)
$$
\sum_{\{i_{1}, \dots, i_{k+1}\}\cup \{j_{1}, \dots, j_{k+1}\}=\{1\sp{m_1}, \dots, \ell\sp{m_\ell}\}} C_{\bf i, j}\,X_{i_{1},j_{1}}(z)\dots X_{i_{k+1},j_{k+1}}(z)=0,
$$
with $C_{\bf i, j}\neq0$, the sum runs over all such partitions of the multiset $\{1\sp{m_{1}}, \dots, \ell\sp{m_{\ell}}\}$. The leading terms of these relations are most conveniently described as monomials  (see \cite{BPT}, Proposition 2)
 \begin{equation}\label{E: leading terms for FS Cell}
 X_{i_t j_t}(-j-1)^{b_{i_t j_t}}\dots  X_{i_1 j_1}(-j-1)^{b_{i_1 j_1}}
X_{i_s j_s}(-j)^{a_{i_s j_s}}\dots X_{i_{t+1} j_{t+1}}(-j)^{a_{i_{t+1} j_{t+1}}},
 \end{equation}
where 
$(i_p,j_p)\neq(i_{p+1}, j_{p+1})$ and
 \begin{equation}\label{E: diagonal paths in Gamma}
i_1 \leq \dots \leq i_t   \leq j_t\leq\dots\leq j_1  \leq i_{t+1} \leq \dots \leq i_s \leq j_{s} \leq \ldots\leq j_{t+1},
 \end{equation}
$$
b_{i_1 j_1}+\ldots+b_{i_t
j_t}+a_{i_{t+1} j_{t+1}}+\ldots+a_{i_s j_s}=k+1.
$$
The factors $X_{pq}$ of leading terms (\ref{E: leading terms for FS Cell}) lie on diagonal paths in $B_1$ with $i=i_{t+1}$:
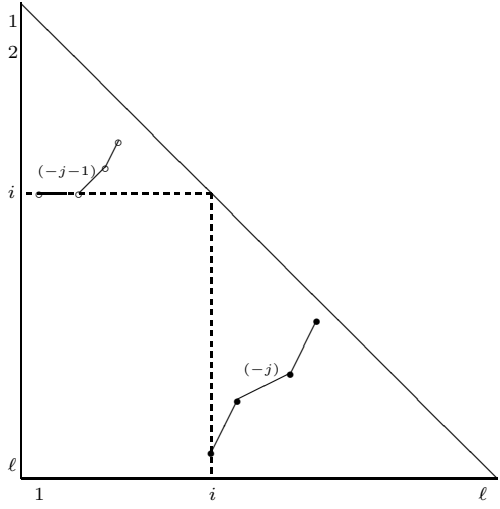
\begin{center}\begin{picture}(200,200)(-10,-10) 
\put(0,0){\line(1,0){180}} \put(0,0){\line(0,1){180}}
\put(180,0){\line(-1,1){180}} 
\put(-5,171){$\scriptstyle 1$}
\put(-5,159){$\scriptstyle 2$}
\put(-5,3){$\scriptstyle \ell$} \put(5,-8){$\scriptstyle
1$}
\put(173,-8){$\scriptstyle \ell$}

\put(84,40){$\scriptscriptstyle (-j)$}

\put(110,58){$\scriptscriptstyle \bullet$}
 \put(102,40){\line(1,2){10}}
\put(100,38){$\scriptscriptstyle \bullet$}
 \put(82,30){\line(2,1){20}}
\put(80,28){$\scriptscriptstyle \bullet$}
 \put(72,10){\line(1,2){10}} 
\put(70,8){$\scriptscriptstyle  \bullet$}

\multiput(72,108)(-4,0){18}{\line(-1,0){2}}
\multiput(72,108)(0,-4){27}{\line(0,-1){2}}

\put(-5,106){$\scriptstyle i$}
\put(71,-8){$\scriptstyle i$}

\put(35,126){$\scriptscriptstyle \circ$}
\put(30,116){$\scriptscriptstyle \circ$}
\put(20,106){$\scriptscriptstyle \circ$}
\put(5,106){$\scriptscriptstyle \circ$}
\put(7,108){\line(1,0){10}} 
\put(32,118){\line(1,2){5}} 
\put(22,108){\line(1,1){10}}
\put(6,115){$\scriptscriptstyle
(-j-1)$}
\end{picture}\end{center}
Note that there is $\ell$ choices for the upper triangle in $B_1$, from $i=1$, when the upper triangle  is $\{11\}$, all the way to $i=\ell$, when the upper triangle is the entire $B_1$. For example, for $\ell=2$ the triangle $B_1$ is
$$\begin{array}{cc}
11 &   \\
12 & 22   .
\end{array}
$$
Then for given $j$ and $i=1$ we have two possible diagonal paths in the lower triangle, vertical and horizontal, with leading terms
\begin{equation}\label{E: leading terms FS C2 1}
\begin{aligned}
&X_{{11}}(-j-1)\sp{b_{11}}X_{{12}}(-j)\sp{a_{12}}X_{{11}}(-j)\sp{a_{11}}\,,
\quad \ \,  \qquad b_{11}+a_{12}+a_{11}= k+1,\\
&X_{{11}}(-j-1)\sp{b_{11}}X_{22}(-j)\sp{a_{{22}}}X_{{12}}(-j)\sp{a_{12}}\,,
\quad \ \, \qquad b_{{11}}+a_{{22}}+a_{{12}}= k+1,\\
\end{aligned}
\end{equation}
while for $i=2$ we have two possible diagonal paths in the upper triangle, vertical and horizontal, with leading terms
\begin{equation}\label{E: leading terms FS C2 2}
\begin{aligned}
&X_{{12}}(-j-1)\sp{b_{{12}}}X_{{11}}(-j-1)\sp{b_{{11}}}X_{{22}}(-j)\sp{a_{22}}\,,
\quad b_{{12}}+b_{{11}}+a_{{22}}= k+1,\\
&X_{{22}}(-j-1)\sp{b_{22}}X_{{12}}(-j-1)\sp{b_{12}}X_{{22}}(-j)\sp{a_{22}}\,,
\quad b_{22}+b_{12}+a_{22}= k+1.
\end{aligned}
\end{equation}
These leading terms are precisely the leading terms (\ref{E: leading terms FS B2}) which led to difference conditions (\ref{E: B2 FS difference conditions}) as a system of inequalities for all $j\in \mathbb N$. In a similar way for the monomial 
$$
\prod X_{pq}(-j)\sp{m_{pq;j}}
$$
we can write difference conditions as a system of inequalities
 \begin{equation}\label{E: FS difference conditions Cell}
 m_{i_1 j_1;j+1}+\ldots+m_{i_t
j_t;j+1}+m_{i_{t+1} j_{t+1}; j}+\ldots+m_{i_s j_s; j}\leq k
 \end{equation}
for all $j\in \mathbb N$ and diagonal paths (\ref{E: diagonal paths in Gamma}) in $B_1$. It is proved in \cite{BPT} that the spanning set of vectors (\ref{E: basis of FS subspace}) satisfying difference conditions is a basis of $W_{C^{(1)}_{\ell}}(k\Lambda_0)$ for all $k\geq 1$ and $\ell\geq2$.

\subsection{Combinatorial coincidence of  $\ell t\left( \bar{R}\sp0_k \right)$ for $C_{2\ell}$ and $\ell t\left( \bar{R}_k \right)$ for $C_{\ell}$}
For $\ell=1$ we have the combinatorial coincidence of  $\ell t\left( \bar{R}\sp0_k \right)$ for $C_{2}\cong B_{2}$ and $\ell t\left( \bar{R}_k \right)$ for $C_{1}\cong A_1$ if we use the identification (\ref{E: identification of bases for B2 and A1}), that is, if we identify two bases
$$
\begin{array}{cc}
11 &  \\
12 & 22 
\end{array}
 \longleftrightarrow \ \
 \begin{array}{cc}
 11 &   \\
{1}\underline{1}& \underline{1}\underline{1}.  
\end{array}
$$
of $A_1$-modules. For $\ell\geq2$ we can identify bases of two vector spaces for finite dimensional Lie algebra representations,
$$
 L_{A_{2\ell-1}}(2\omega_1)\qquad \text{and} \qquad   L_{C_\ell}(\theta),
$$
the basis (\ref{E: basis of g1}) for $C_{2\ell}$ and the basis  (\ref{E: basis of g}) for $C_{\ell}$. For example, for $\ell=3$ we identify two bases
\begin{equation}\label{E: identification of bases for ell geq 2}
\begin{array}{cccccc}
11 &  &&  & & \\
12 & 22 & & & & \\
13 & 23 & 33 & & & \\
14 & 24 & 34 &44 & & \\
15 & 25 & 35 &45 &55 & \\
16 & 26 & 36 & 46&56 & 66
\end{array}
 \longleftrightarrow \ 
 \begin{array}{cccccc}
11 &  &&  & & \\
12 & 22 & & & & \\
13 & 23 & 33 & & & \\
1\underline{3} & 2\underline{3} & 3\underline{3} & \underline{3}\underline{3}  & & \\
1\underline{2} & 2\underline{2} & 3\underline{2} & \underline{3}\underline{2}  & \underline{2}\underline{2}& \\
1\underline{1} & 2\underline{1} & 3\underline{1} &
\underline{3}\underline{1}  & \underline{2}\underline{1} &
\underline{1}\underline{1}.
\end{array}
\end{equation}
We consider (cf. \cite{Bo})
$$
\mathfrak{sp}_{2\ell}\subset\mathfrak{sl}_{2\ell}.
$$
Then the Weyl dimension formula gives the following:
\begin{proposition}\label{proposition}
For all $m\in\mathbb N$ the restriction of the representation $L_{A_{2\ell-1}}(2m\omega_1)$ for the Lie algebra $\mathfrak{sl}_{2\ell}$ to the subalgebra $\mathfrak{sp}_{2\ell}$ remains irreducible, i.e., as a vector space
$$
L_{C_\ell}(m\theta)=L_{A_{2\ell-1}}(2m\omega_1)\vert_{C_\ell}\cong S\sp{2m}(\mathbb C\sp{2\ell})\vert_{C_\ell}.
$$
In particular, under the identification (\ref{E: identification of bases for ell geq 2}) of two bases and the respective orders, the Feigin-Stoyanovsky's type subspace $W_{C^{(1)}_{2\ell}}(k\Lambda_0)$  and the standard module $L_{C^{(1)}_{\ell}}(k\Lambda_0)$ have the same combinatorial description of the  leading terms of relations for all $k\geq1$ and $\ell\geq1$.
\end{proposition}

\begin{remark}  The above argument gives an alternative construction of the set of leading terms of relations for standard modules $\ell t\left( \bar{R}_k \right)$ for $C_{\ell}$, first obtained combinatorialy in \cite{PS2}.
\end{remark}

\begin{remark} So far it seems that it is easier to prove linear independence of spanning sets (\ref{E: basis of FS subspace}) for Feigin-Stoyanovsky's type subspaces than to prove linear independence of spanning sets (\ref{E: basis of standard module}) for standard modules. So the question is: can one use the results for $W_{C^{(1)}_{2\ell}}(k\Lambda_0)$ to study $L_{C^{(1)}_{\ell}}(k\Lambda_0)$ for $k\geq2$\,? 
\end{remark}

\begin{remark} D. Adamovi\' c pointed out that the isomorphism $L_{A_{2\ell-1}}(n\omega_1)\vert_{C_\ell} \cong L_{C_\ell}(n\omega_1)$ appears in \cite{AP} in the context of vertex algebras (cf. Proposition 7) and may be viewed as a consequence of the quantum Galois theory \cite{DM}.
\end{remark}

\end{document}